\documentstyle{amsppt}

\magnification=1100 \NoBlackBoxes\nologo 
\vsize=19.5cm

\def\inv{^{-1}}
\def\?{{\bf{??}}}

\def\Hilb{\text{Hilb}}
\def\fHilb{\text{fHilb}}

\def\A{\Bbb A}

\def\C{\Bbb C}

\def\P{\Bbb P}

\def\Spec{\text{\rm Spec} }

\def\[{\big[}
\def\]{\big]}

\def\O{\Cal O}

\def\Sym{\text{Sym}}

\def\m{\frak m}
\def\k{\frak k}

\def\1/2{\frac{1}{2}}

\def\I{\Cal I}

\def\2{{[2]}}

\topmatter
\title A note on Hilbert schemes of Nodal Curves
\endtitle
\author
Ziv Ran
\endauthor

\date 20040215\enddate

\address University of California, Riverside\endaddress
\email ziv\@math.ucr.edu\endemail
\rightheadtext { Hilbert scheme of Nodal Curves} \leftheadtext{Z.
Ran} \abstract We study the Hilbert scheme and punctual Hilbert
scheme of a nodal curve, and the relative Hilbert scheme of a
family of curves acquiring a node. The results are then extended
to flag Hilbert schemes, parametrizing chains of subschemes. We
find, notably, that if the total space $X$ of a family $X/B$ is
smooth (over an algebraically closed field $\k$), then the
relative Hilbert scheme $Hilb_m(X/B)$ is smooth over $\k$ and the
flag Hilbert schemes are normal and locally complete intersection,
but generally singular .
\endabstract
 \thanks \raggedright {Updates and corrections
available at $\underline{math.ucr.edu/\tilde{\ }
ziv/papers/hilb.pdf}$}\linebreak Research Partially supported by
NSA Grant MDA904-02-1-0094; reproduction and distribution of
reprints by US government permitted.
\endthanks

\endtopmatter\document
The Hilbert scheme parametrizes ideal sheaves or subschemes of
Projective Space or more generally, or a fixed scheme $X$. Perhaps
the simplest case is where $X$ is a smooth curve, for then the
Hilbert scheme $\Hilb_m(X)$ parametrizing length-$m$ subschemes of
$X$ coincides with the symmetric product $\Sym^m(X)$. Our purpose
in this note is to study what is in a sense the next simplest
case, where $X$ is essentially a curve with ordinary nodes, with
planar equation formally equivalent to $xy=0.$ Besides
$\Hilb_m(X)$ itself, there are (at least) 2 other types of Hilbert
scheme of natural interest here: the {\it punctual} one
$\Hilb^0_m(X)$, parametrizing length-$m$ subschemes supported at a
node; and the {\it relative} one $\Hilb_m(\tilde{X}/B)$,
parametrizing length-$m$ subschemes in the fibres of the family
$\tilde{X}/B$ that is the map $\A^2\to \A^1$ given by $xy=t$ (both
$\Hilb_m(X)$ and $\Hilb_m(\tilde{X}/B)$ may conveniently be viewed
as germs, or formal schemes, along $\Hilb^0_m(X)$. All these will
be determined below. We find that $\Hilb^0_m(X)$ is a connected
chain of $(m-1)$ nonsingular rational curves meeting
transeversely; $\Hilb_m(X)$ is a connected chain of $(m+1)$
nonsingular $m$-dimensional germs, whose first and last members
are supported on points and each other member is supported on a
component of $\Hilb^0_m(X)$, and only consecutive members
intersect (and those transversely); finally,
$\Hilb_m(\tilde{X}/B)$ is a smooth $(m+1)$-fold.\par To elucidate
the relations between the $\Hilb_m$ for different $m$, with an eye
to enumerative applications based on recursion on $m$, we will
also study some flag Hilbert schemes $\Hilb_{m.}$, parametrizing
chains of ideals whose colengths form a given sequence
$m.=(m_1>m_2>...)$. For some purposes, these are easier to work
with than ordinary Hilbert schemes, due to the natural maps
between the $\Hilb_{m.}$ for different $(m.)$. For example, for
$m.=(m,m-1)$, we find that the punctual Hilbert scheme
$\Hilb^0_{m.}(X)$ is a chain of $2m-3$ copies of $\P^1$ that
alternate between those coming from $\Hilb^0_m(X)$ and from
$\Hilb^0_{m-1}(X)$. The relative Hilbert scheme
$\Hilb_{m.}(\tilde{X}/B)$ is still smooth. Things begin to change
though with $m.=(m,m-1,m-2)$. Here $\Hilb^0_{m.}(X)$ is still a
chain of $2m-3$ components, but now only the 2 external ones on
each side are $\P^1$'s and the rest are $\P^1\times\P^1$. The
relative Hilbert scheme $\Hilb_{m.}(\tilde{X}/B)$ is no longer
smooth, but it is still normal with locally complete intersection
singularities. A similar picture emerges for longer chains,and in
particular for the 'full flag' case
$\fHilb_m=\Hilb_{m,m-1,...,1}.$\par In \cite{R}, which uses some
of these results, we will develop an alternative and more
'geometric' approach to these objects. In particular, we will
identify $\fHilb_m(\tilde{X}/B)$ with a certain space
$W^m(\tilde{X}/B)$ constructed as an explicit blow-up of the
relative Cartesian product $\tilde{X}^m/B$. The proof of this
identification uses some of the results here. We will then apply
the spaces $W^m(\tilde{X}/B)$ to some enumerative problems. In
these applications, the relatively simple relationship between
$\fHilb_m$ and $\fHilb_{m-1}$ is critical.\par \comment
 In this section we will first study the
punctual Hilbert scheme of length-$r$ schemes supported at the
origin on a germ of a node $xy=0.$ Then we will study the full
Hilbert scheme of this germ. Finally we will study the relative
Hilbert scheme of a flat family of nodal curves in the
neighborhood a node. These results will form the local foundation
on which we shall in \S 2 construct our global parameter spaces
$W^r(X/B)$.\par
\endcomment
We work over a fixed algebraically closed field $\k$. We denote by
$R$ localizarion of the ring
$$\k[x,y]/(xy)$$ and at its maximal ideal
$(x,y)$. A typical element of $R$ can be written in the form
$$u({a+\sum\limits_{i\geq 1}b_ix^i+c_iy^i})$$ where $u$ is a unit
and the sum is finite.
  The formal completion
$$\hat{R}=\k[[x,y]]/(xy)=\{{a+\sum\limits_{i\geq 1}b_ix^i+c_iy^i}\}$$
(sum not necessarily finite) is isomorphic to the formal
completion of the local ring at any 1-dimensional ordinary node.
We seek first to determine the punctual Hilbert scheme
Hilb$^0_m(R)$ of colength-$m$ ideals in $R$ which, as is well
known, is naturally isomorphic to Hilb$^0_m(\hat{R})$. At this
point, we do not seek to define or compute a natural scheme
structure on $\Hilb^0_m(R)$ (so calling it the punctual Hilbert
{\it scheme} is something of a misnomer)-- this will be done later
(see Corollary 8). For now we simply view $\Hilb^0_m(R)$ as an
algebraic set endowed with a flat family of ideals that yields a
bijective correspondence between the (closed, $\k-$valued) points
of $\Hilb^0_m(R)$ and the colength-$m$ ideals of $R$. \proclaim
{Theorem 1} (i) Every ideal $I<R$ of colength $m$ is of one of the
following, said to be of type ${(c^m_i)}, (q^m_i)$, respectively:
$$ I^m_i(a)=
(y^{i}+ax^{m-i}), 0\neq a\in\k , i=1,...,m-1;$$
$$Q^m_i=(x^{m-i+1},y^{i}),
i=1,...,m.$$ (ii) The closure $C^m_i$ in the Hilbert scheme of the
set of ideals of type $(c^m_i)$ is isomorphic to $\P^1$ and
consists of the ideals of types $(c^m_i)$ or $(q^m_i)$ or
$(q^m_{i+1})$. In fact, we have
$$\lim\limits_{a\to 0}I^m_i(a)=Q^m_i,$$
$$\lim\limits_{a\to \infty}I^m_i(a)=Q^m_{i+1}.$$
(iii) The punctual Hilbert scheme Hilb$_m^0(R)$, as algebraic set,
is a rational chain
$$C^m_1\cup_{Q^m_2}C^m_2\cup\cdots\cup_{Q^m_{m-1}} C^m_{m-1};\tag 1$$
it has ordinary nodes at $Q^m_2,...,Q^m_{m-1}$ and is smooth
elsewhere. \comment \linebreak (iv) The only colength-$(m-1)$
ideal containing $I^m_i(a)$ is $Q^{m-1}_i$; the only
colength-$(m-1)$ ideals containing $Q^m_i$ are the
$I^{m-1}_{i-1}(a)$ for $a\neq 0$ and their limits $Q^{m-1}_i$ and
$Q^{m-1}_{i-1}$.
\endcomment

\endproclaim
\demo{proof} Recall that elements $z, z'\in R$ are said to be {\it
associate} if $z=uz'$ for some unit $u$. Note that any nonzero
nonunit $z\in R$ is associate to a uniquely determined element of
the form $x^\alpha$ or $y^\beta$ or $x^\alpha+ay^\beta, a\neq 0,
\alpha,\beta>0,$ in which case we will say that $z$ is of type
$(\alpha,0)$ or $(0,\beta)$ or $(\alpha,\beta)$, respectively.
Note also that for any ideal $I$ of colength $m$ we have
$$x^m, y^m\in I.$$
Now given $I$ of colength $m$, pick $z\in I$ of minimal type
$(\alpha,\beta)$, with respect to the natural partial ordering on
types. Suppose to begin with that $\alpha,\beta>0.$ Then note that
$$x^{\alpha+1}, y^{\beta+1}\in I,$$
and consequently $(\alpha,\beta)$ is unique: indeed if
$(\alpha',\beta')$ is also minimal then we may assume
$\alpha'>\alpha$, hence $x^{\alpha'}\in I$, hence $y^{\beta'}\in
I,$ contradicting minimality. Hence $(\alpha,\beta)$ is unique and
it is then easy to see that the element $z'=x^\alpha+ay^\beta\in
I$ is unique as well, so clearly $z'$ generates $I$ and $I$ is of
type $(c^m_\beta)$.\par Thus we may assume that any minimal
element of $I$ is of type $(\alpha,0)$ or $(0,\beta)$. Since
$x^m,y^m\in I$, $I$ clearly contains minimal elements of type
$(\alpha,0)$ and $(0,\beta)$, and then it is easy to see that $I$
is of type $(q^m_\beta)$. This proves assertion (i). Since
$I^m_i(a)$ contains $y^{i+1}, x^{m-i+1},$ assertion (ii) is easy.
As for (iii), let $C=\bigcup\limits_iC^m_i$ be an abstract nodal
curve as in (1). It follows from (ii) that each $C^m_i$ carries a
flat family of ideals, i.e. admits a natural morphism- clearly
injective- to $\Hilb_m(R)$. Since these morphisms agree on the
intersections $C^m_i\cap C^m_{i+1}=Q^m_{i+1},$ they yield a
morphism, again clearly injective, from $C$ to $\Hilb_m(R)$, which
identifies $C$ with $\Hilb^0_m(R)$ as an algebraic
set.\qed\enddemo \comment
\par For (iii), we note that the punctual Hilbert scheme
parametrizes punctual deformations of an ideal $I<R$, parametrized
by a local augmented artin $\k-$algebra $S$. This means ideals
$$I_S<R_S:=R\otimes S$$
such that $I_S\subseteq (x,y)$, $R_S/I_S$ is $S$-free of rank $m$,
and which are punctual in the sense that $I_S$ is contained in a
unique ideal $J$ that is maximal subject to
$$J\cap S=0$$ (in which
case $J=(x,y)$ automatically). The assertion is that this scheme
is 1-dimensional and is smooth at at any $I^m_i(a), i=1,...,m-1,
Q_1^m$ and $Q_m^m$ and nodal at $Q^m_i, 1<i<m.$
 We consider the case of $Q^m_i, 1<i<m$
as the others are similar but simpler. Then $I_S$ is generated by
$$f=x^{m+1-i}+f_1(x)+f_2(y),$$
$$g=y^i+g_1(x)+g_2(y)$$
where the $f_i,g_j$ are polynomials (or power series) without
constant term (by punctuality) and with coefficients in $\m_S$.
Denote the lowest term of $f_1$ by $bx^j, j>0.$ Suppose $
j<m+1-i.$ Then note that
$$I_S\subseteq (x^{m+1-i-j}+b,y, g_1(x)),$$
and the latter is a proper ideal by $ j<m+1-i,$ is not contained
in $(x,y)$ if $b\neq 0,$ and  meets $S$ trivially because $g_1$
has coefficients in $\m_S$. This contradicts the punctuality of
the deformation $I_S$ (i.e. the uniqueness of $J=(x,y)$ as above).
Therefore $j\geq m+1-i$ and after multiplying $f$ by a unit we may
assume
$$f=x^{m+1-i}+by^r.$$ Similarly, we may assume
$$g=y^i+cx^s.$$
Note
if $r\geq i$ then replacing $f$ by $f-bx^{r-i}g$ yields a
polynomial in $x$ only which by the above remarks may be assumed
to be $x^{m+1-i}$, so we can take $b=0$; hence we may assume $r<i$
and similarly $s<m+1-i$. On the other hand if $r<i-1$ then
$by^{i-1}\in I_S$ whereas by flatness,
$$1, x,...,x^{m-i}, y,...,y^{i-1}$$
must be an $S$- free basis of $R_S/I_S$. Therefore $r=i-1, s=m-i$,
i.e. $$I_S=(x^{m+1-i}+by^{i-1}, y^i+cx^{m-i}).\tag 2$$ Finally,
multiplying $f,g$ by $c, b, y, x$ we get
$$cx^{m+1-i}+bcy^{i-1}, by^i+bcx^{m-i}, by^i, cx^{m-i+1}\in I_S,$$
hence $bcy^{i-1}, bcx^{m-i}\in I_S,$ so flatness implies $bc=0.$
Thus, punctual $S-$deformations of $Q^m_i$ are parametrized by a
pair $b,c,\in\m_S$ with
$$bc=0,\tag 3$$ i.e. by the local $\k-$algebra homomorphisms of the
local ring of a node to $S$, so the punctual Hilbert scheme itself
is a 1-dimensional (ordinary) node locally at $Q^m_i$, as claimed.
\qed\enddemo
\endcomment

Next we determine the structure of the full Hilbert scheme of $R$
and $\hat{R}$: \proclaim{Theorem 2} The Hilbert scheme
$\text{Hilb}_m(R)$ (resp. $\text{Hilb}_m(\hat{R})$ is a chain
$$D^m_0\cup D^m_1\cdots D^m_{m-1}\cup D^m_m$$
where each $D^m_i$ is a smooth and $m-$dimensional germ (resp.
formal scheme) supported on $C^m_i$ for $i=1,...,m-1$ or $Q^m_i$
for $i=0,m$; for $i=1,...,m-1, D^m_i$ meets its neighbors
$D^m_{i\pm 1}$ transversely in dimension $m-1$ and meets no other
$D^m_i.$ The generic point of $D^m_i$ corresponds to subscheme of
$\Spec(R)$ comprised of $m-i$ points on the $x$-axis and $i$
points on the $y$-axis.\endproclaim \demo{proof} Clearly
Hilb$_m(R)$ (resp. $\Hilb_m(\hat{R})$) is a germ (resp. formal
scheme) supported on Hilb$^0_m(R)$, so this is a a matter of
determining the scheme structure of Hilb$_m(R)$ and
$\Hilb_m(\hat{R})$ at each point of Hilb$^0_m(R)$, which may be
done formally by testing on Artin local algebras. Again, we shall
do so at $Q^m_i, i>1$ as the cases of $I^m_i(a)$ and $Q^m_1$ are
similar and simpler. Given $S$ artinian local augmented, a  flat
$S$-deformation of $I=Q^m_i$ is given by an ideal
$$I_S=(f,g),$$
$$f=x^{m+1-i}+f_1(x)+f_2(y),\tag 4$$
$$g=y^i+g_1(x)+g_2(y),$$
where $f_i,g_j$ have coefficients in $\m_S$, and such that
$R_S/I_S$ is $S-$ free of rank $m$, in which case it is clear by
Nakayama's Lemma that
$$1, x,..., x^{m-i}, y ,..., y^{i-1}$$
is a free basis for $R_S/I_S$. It is easy to see that we may
assume $f_1, g_1$ are in fact polynomials of degree $\leq m-i$ and
$f_2,g_2$ are of degree $< i$ and $f_2,g_2$ have no constant term.
Let's write
$$f_1(x)=\sum\limits_0^{m-i} a_jx^j, f_2(y)=\sum\limits
_1^{i-1} b_jy^j,\tag 5$$
$$g_1(x)=\sum\limits_0^{m-i} c_jx^j, g_2(y)=\sum\limits
_1^{i-1} d_jy^j.\tag 6$$ Now obviously
$$yf-b_{i-1}g\equiv 0\equiv xg-c_{m-i}f\mod I_S.$$ Writing
these elements out in terms of $1,x,...,x^{m-i}, y,...,y^{i-1}$
yields relations among $1,x,...,x^{m-i}, y,...,y^{i-1}$. Since the
latter elements form an $S$-free basis of $R_S/I_S$, those
relations must be trivial. In other words, we have exact
equalities rather than congruences:
$$yf-b_{i-1}g= 0= xg-c_{m-i}f.$$
Writing out this equality term by term yields the following
identities
$$b_j=b_{i-1}d_{j+1}, j=1,...,i-2,$$
$$b_{i-1}d_1=a_0,$$
$$b_{i-1}c_j=0, j=0,...,m-i,\tag 7$$
$$c_j=c_{m-i}a_{j+1}, j=0,...,m-i-1,$$
$$c_{m-i}a_0=0,$$
$$c_{m-i}b_j=0, j=1,...,i-1.$$
(If $i=1$ only lines 4 and 5 of display (7) are operational.)
Conversely, suppose the relations (7) are satisfied, or
equivalently
$$yf-b_{i-1}g= 0= xg-c_{m-i}f.\tag 8$$
By Nakayama's Lemma, $1, x,...,x^{m-i}, y,...,y^{i-1}$ generate
$R_S/I_S$, hence to show (4) defines a flat family it suffices to
show these elements admit no nontrivial $S$-relations mod $I_S$.
To this end,, suppose
$$u_{m-i}(x)+v_{i-1}(y)=A(x,y)f+B(x,y)g\tag 9$$
where $u, v, A, B$ are all polynomials in the indicated variables
and of the indicated degees (if any) with coefficients in $S$ and
$v$ has no constant term; in fact it clear a priori that then
$A,B$ must have coefficients in $\m_S$. Then the relations (8)
 allow us
to rewrite (9) as
$$u_{m-i}(x)+v_{i-1}(y)=A'(x)f+B'(y)g,\tag 9'$$
with $A'(x)=\sum a'_kx^k\in\m_S[[x]], B'(x)=\sum b'_ky^k\in
\m_S[[y]].$ Comparing coefficients of $x^{m-i+1},...,y^{i},...$ in
(9') we get relations
$$-a'_0=a'_1a_{m-i}+..., -b'_0=b'_1b_{i-1}+...$$
$$...\tag 10$$
$$-a'_k=a'_{k+1}a_{m-i}+..., -b'_k=b'_{k+1}b_{i-1}+...$$
Starting from the fact that $a_j, a'_k\in\m_S, \forall j,k$, we
infer from these first that, in fact, $a'_k\in\m_S^2, \forall k$;
plugging the latter fact back into the relations (10) we then
infer $a'_k\in\m_S^3, \forall k,$ and so on. Since $S$ is artinian
it follows that $a'_k=0, \forall k$ and likewise for $b'_k.$ Thus
$$A'=B'=u_{m-i}=v_{i-1}=0,$$
hence there are no nontrivial relations, as claimed.\par Thus the
Hilbert scheme is embedded in the space of the variables
$$a_1,...,a_{m-i}, d_1,...,d_{i-1}, b_{i-1},c_{m-i},$$ i.e.
$\A^{m+1}$, and defined by the relation
$$b_{i-1}c_{m-i}=0.\tag 11$$
Thus it is a union of 2 smooth $m-$dimensional components meeting
transversely in a smooth $(m-1)-$dimensional subvariety. The
generic point on the component where $b_{i-1}=0$ (resp.
$c_{m-i}=0$) is clearly an ideal generated by $g$ (resp. $f$),
which has the properties as claimed.
\par
In the case $I=I^m_i(a)=(y^i+ax^{m-i})$, a similar analysis shows
that an $S$-deformation of $I$ is given by the principal ideal
$$I_S=(y^i+\tilde{a}x^{m-i}+f_1(x)+f_2(y))$$
where
$$\tilde{a}\in S, \tilde{a}\equiv a\mod\m_S,$$
$$f_1(x)=\sum\limits_0^{m-i-1} a_jx^j, f_2(y)=\sum\limits
_1^{i-1} b_jy^j,$$
$$a_j, b_j\in\m_S,$$
and via $(\tilde{a}, a_0,...,a_{m-i-1}, b_1,...,b_{i-1})$ we may
identify Hilbert scheme locally with $\A^m.$ \qed
\enddemo
\remark{Remark 2.1} We note that in terms of the above coordinates
the subset $\Hilb^0_m(R)\subset \Hilb_m(R)$ is defined by
$$a_1=...=a_{m-i}=d_1=...=d_{i-1}=0,$$i.e. by the conditions
$$f(x,0)=x^{m-i+1}, g(0,y)=y^i.$$ We shall see below that this, in
fact, defined the 'natural' scheme structure on
$\Hilb^0_m(R)$.\qed
\endremark
 Next we consider the relative local
situation, i.e. that of a germ of a (1-parameter) family of curves
with smooth total space specializing to a node. Thus set
$$\tilde{R}= \k[x,y]_{(x,y)},    B=\k[t]_{(t)},
$$
and view $\tilde{R}$ as a $B-$module via $xy=t$. As is well known,
this is the versal deformation of the node singularity $xy=0$, so
any family of nodal curves is locally a pullback of this.
\proclaim{Theorem 3} The relative Hilbert scheme
$\text{Hilb}^m(\tilde{R}/B)$ is formally smooth, formally
$(m+1)-$dimensional over $\k$.\endproclaim \demo{proof} The
relative Hilbert scheme parametrizes length-$m$ schemes contained
in fibres of Spec$\tilde{(R)}\to$ Spec$(B)$. This means ideals
$I_S<\tilde{R}_S$ of colength $m$ containing $xy-s$ for some
$s\in\m_S$, such that $\tilde{R}_S/I_S$ is $S-$free. The analysis
of these is virtually identical to that contained in the proof of
Theorem 2, except that the relation $b_{i-1}c_{m-i}=0$ gets
replaced by
$$b_{i-1}c_{m-i}=s\tag *$$ and lines 3,5,6 of display (7) are replaced,
respectively, by $$b_{i-1}c_j=sa_{j+1}, j=0,...,m-i-1$$
$$c_{m-i}a_0=sd_1$$$$c_{m-i}b_j=sd_{j+1}, 1\leq j\leq i-2,$$
relations which already follow from the other relations (in lines
1,2,4 of display (7)) combined with (*).
 Thus, the relative Hilbert scheme
is the subscheme of the affine space of the variables
$a_1,...,a_{m-i},d_1,...,d_{i-1}, b_{i-1}, c_{m-i}, t$ defined by
the relation
$$b_{i-1}c_{m-i}=t,\tag 12$$
hence is smooth as claimed.\qed\enddemo \remark{Remark 3.1}After
this was written, the author was informed by Prof. I. Smith of
some related work by himself and Prof. S. Donaldson \cite{DS, Sm}
which considers the relative Hilbert scheme of a pencil of nodal
curves on a smooth surface from a rather different viewpoint,
valid in the symplectic category over $\C$; in particular, they
prove in this context an analogue of Theorem 3 (smoothness of the
total space of the relative Hilbert scheme).\endremark
\remark{Construction 3.2}An explicit construction of $\Hilb_m(R)$,
globally along $\Hilb^0_m(R)$, can be given as follows. Let
$C_1,...,C_{m-1}$ be copies of $\P^1$, with homogenous coordinates
$u_i,v_i$ on the $i$-th copy. Let $\tilde{C}\subset
C_1\times...\times C_{m-1}\times\A^1$ be the subscheme defined by
$$v_1u_2=tu_1v_2,..., v_{m-2}u_{m-1}=tu_{m-2}v_{m-1}.$$ Thus
$\tilde{C}$ is a reduced complete intersection of divisors of type
$(1,1,0,...,0), (0,1,1,0,...,0)$ ,..., $(0,...,0,1,1)$ and it is
easy to check that the fibre of $\tilde{C}$ over $0\in\A^1$ is
$$\tilde{C}_0=\bigcup\limits_i [1,0]\times...\times[1,0]\times
C_i\times[0,1]\times...\times[0,1]$$ and that in a neighborhood of
$\tilde{C}_0$, $\tilde{C}$ is smooth and $\tilde{C}_0$ is its
unique singular fibre over $\A^1.$ We may identify $\tilde{C}_0$
is an obvious way with $\Hilb_m^0(R).$ Next consider an affine
space $\A^{2m}$ with coordinates $a_0,...,a_{m-1},
d_0,...,d_{m-1}$ and let $\tilde{H}\subset\tilde{C}\times\A^{2m}$
be the subscheme defined by
$$a_0u_1=tv_1, d_0v_{m-1}=tu_{m-1}$$
$$a_1u_1=d_{m-1}v_{1},...,a_{m-1}u_{m-1}=d_1v_{m-1}.$$
Set $L_i=p_{C_i}^*\O(1).$ Then consider the subscheme of
$\tilde{H}\times_{\A^1}\Spec(\tilde{R})$ defined by the equations
$$F_0:=x^m+a_{m-1}x^{m-1}+...+a_1x+a_0\in
\Gamma(\O_{\tilde{H}}\otimes R)$$
$$F_1:=u_1x^{m-1}+u_1a_{m-1}x^{m-2}+...+u_1a_2x+u_1a_1+v_1y
\in\Gamma(L_1\otimes R)$$ ...
$$F_i:=u_ix^{m-i}+u_ia_{m-1}x^{m-i-1}+...+u_ia_{i+1}x+u_ia_i+
v_id_{m-i+1}y+...+v_id_{m-1}y^{i-1}+ v_iy^i$$
$$
\in\Gamma(L_i\otimes R)$$ ...
$$F_m:=d_0+d_1y_1+...+d_{m-1}y^{m-1}+y^m\in
\Gamma(\O_{\tilde{H}}\otimes R).$$ In view of the above results,
it is easy to check that the ideal sheaf $\I$ generated by
$F_0,...,F_m$ defines a subscheme of $\tilde{H}\times\Spec(R)$
that is flat over $\tilde{H}$ and that may serve to identify
$\tilde{H}$ with $\Hilb(\tilde{R})/B.$ Locally at a point of type
$c^m_i$ (resp. $Q^m_i), \I$ is generated by $F_i$ (resp. $F_{i-1},
F_i$).\par Now let $E$ be the universal bundle on $\tilde{H}$,
whose fibre at a point corresponding to a length-$m$ scheme $z$ is
$\Gamma(\O_z),$ and let $V$ be the trivial rank$-(2m+1)$ bundle on
the symbols $1, x,...,x^m,y,...,y^m.$ Since  $1,
x,...,x^m,y,...,y^m$ generate $\Gamma(\O_z)$ for all
$z\in\Hilb_m(\tilde{R}),$ we have a surjection $V\to E$. Then
$F_0,...,F_m$ together yield a map
$$F:\O_{\tilde{H}}\oplus\bigoplus\limits_1^{m-1}L_i\inv\oplus\O_{\tilde{H}}\to
V$$so we get an exact sequence
$$0\to\O_{\tilde{H}}\oplus\bigoplus\limits_1^{m-1}L_i\inv\oplus\O_{\tilde{H}}
\to V\to E\to 0.$$ In particular, it follows that
$c_1(E)=L_1+...L_{m-1}.$ \par Note that
$\A^m_{a_0,...,a_{m-1}},\A^m_{d_0,...,d_{m-1}}$ may be identified,
respectively, with $\Sym^m\A^1_x,$ $ \Sym^m\A^1_y$ and accordingly
it will be convenient to use the notation $\sigma_i=a_{m-i},
\tau_i=d_{m-i}$, these being respectively the $i$-th elementary
symmetric functions in the roots of $F_0, F_m.$ I claim, in fact,
that the projection $c:\tilde{H}\to \A^{2m+1}$ may be identified
with the cycle map $\tilde{H}\to\Sym^m(\Spec(\tilde{R})/B):$ this
is because the natural map
$$\Phi:\Sym^m(\Spec(\tilde{R})/B)\to\A^m_\sigma\times\A^m_\tau\times B$$
is an embedding. $\Phi$ is an embedding because its restrictions
on the generic fibre and the special fibre over $B$ are
embeddings. Moreover, it is not hard to see that the image of
$\Phi$ (which coincides with that of $c$) is scheme-theoretically
defined by the equations
$$\sigma_i\tau_j=t\sigma_{i-1}\tau_{j-1},\ \ \forall i, j =1,...,m,
\ \  i+j>m$$ where we set $\sigma_0=\tau_0=1.$ Setting
$S_m=$im$(\Phi),$ we have that $\tilde{H}$ is a Weil divisor in
$S_m\times\tilde{C}$.

 \qed\endremark

 The local analysis immediately yields some
conclusions for the Hilbert scheme of a nodal curve:
\proclaim{Corollary 4} Let $C_0$ be a curve with only $k$ nodes as
singularities and $c$ irreducible components. Then\par\noindent
(i)\ \ $\text{Hilb}_m(C_0)$ is reduced and has precisely
$\binom{m+c-1}{m}$ components, the general element of each of
which corresponds to a reduced subscheme of the smooth part of
$C_0;$
\par\noindent (ii)\ \ let $I$
be a point of $\text{Hilb}_m(C_0)$ having colength $m_i$ at the
$i$-th node of $C_0;$ then locally at $I$,$\text{Hilb}_m(C_0)$ is
a cartesian product of $k$ factors, each of which is a 2-component
normal crossing of dimension $m_i, i=1,...,k$, or a point if
$m_i=0$, times a smooth factor.
\par\noindent (iii)\ \ the fibre of the cycle map
$$cyc:\text{Hilb}_m(C_0)\to \Sym^m(C_0)$$ (cf. \cite{A})
over a cycle having multiplicity $m_i$ at the $i$th node is a
product of 1-dimensional rational chains of length $m_i-1.$
\endproclaim
\demo{proof} It is clear from the explicit analysis in the proof
of Theorem 2 that any subscheme of $C_0$ deforms to a reduced
subscheme supported on the smooth part. Such subschemes are
parametrized by an open dense subset of the symmetric product
$\Sym^m(C_0).$ This clearly yields (i) and (ii), while (iii)
follows from the fact that the fibres of $cyc$ are products of
punctual Hilbert schemes.\qed\enddemo Next we extend (most of) the
above results to the case of {\it flag} Hilbert schemes (see
\cite{Se} for a general discussion of those). By definition, for
any decreasing sequence of positive integers
$$ m_.= (m_1>...>m_k),$$ the flag Hilbert scheme Hilb$_{m.}(R)$
parametrizes nested chains of ideals $$I_S^1<...<I_S^k<R_S\tag
13$$ such that $R_S/I^j_S$ is $S$-free of rank $m_j, j=1,...,k$
(i.e. such that each $I_S^j$ is an $S-$ point of Hilb$_{m_j}(R)$);
similarly for punctual and relative Hilbert schemes. Thus a flag
Hilbert scheme is by definition a subscheme of a product (or, in
the relative case, a fibre product) of its 'constituent' ordinary
Hilbert schemes and as such comes equipped with forgetful
projection morphisms to those constituents; moreover the defining
equations for Hilb$_{m.}(R)$ in
$\prod\limits_{j=1}^k$Hilb$_{m_j}(R)$ are just the conditions that
the inclusions (13) hold, and these equations involve only pairs
of successive factors Hilb$_{m_j}(R)$, Hilb$_{m_{j+1}}(R),
j=1,...,k-1$ . In the case of 'full flags', i.e. the case
$$m_.=(m>m-1>...>1),$$ we will denote Hilb$_{m.}(R)$ by
fHilb$_m(R)$. We begin by analyzing the case of pairs of ideals of
relative colength 1: \proclaim{Theorem 5} (i)The only
colength-$(m-1)$ ideal containing $I^m_i(a)$ is $Q^{m-1}_i$; the
only colength-$(m-1)$ ideals containing $Q^m_i$ are the
$I^{m-1}_{i-1}(a)$ for $a\neq 0$ and their limits $Q^{m-1}_i$ and
$Q^{m-1}_{i-1}$.\par\noindent (ii) The punctual flag Hilbert
scheme  Hilb$_{m,m-1}^0(R)$, as algebraic set, is a chain of
nonsingular rational curves of the form
$$C^m_1\cup_{(Q^m_2,Q^{m-1}_1)}C^{m-1}_1\cup_{(Q^m_2,Q^{m-1}_2)}
C^m_2\cup\cdots\cup_{(Q^m_{m-1}, Q^{m-1}_{m-1})} C^m_{m-1};$$ it
has ordinary nodes at $Q^m_2,...,Q^m_{m-1}$ and is smooth
elsewhere. Each component $C^m_i$ projects isomorphically to its
image in Hilb$_m^0(R)$ and to a point $Q^{m-1}_i$ in
Hilb$_{m-1}(R)$, and vice-versa for $C^{m-1}_i.$\par\noindent
(iii) The flag Hilbert scheme $\Hilb_{m,m-1}(R)$, as formal scheme
along  $\Hilb^0_{m,m-1}(R)$, has normal crossing singularities and
at most triple points. Each of its components is formally smooth,
$m-$dimensional and of the form
$$D^{m,m-1}_{i,i'}, i=0,...,m, i-1\leq i'\leq i,$$ and
projects to $D^m_i, D^{m-1}_{i'}$ respectively (cf. Theorem 2);
$D^{m,m-1}_{i,i'}$ meets $D^{m,m-1}_{j,j'}$ nontrivially iff
$$|i-j|+|i'-j'|\leq 1; $$ the components of $\Hilb^0_{m,m-1}(R)$
contained in $D^{m,m-1}_{i,i}$ (resp. $D^{m,m-1}_{i,i-1}$) are
$C^{m}_i$ and $C^{m-1}_i$ (resp. $C^{m-1}_{i-1}$ and $C^{m}_i$).
\par\noindent (iv) The relative flag Hilbert scheme
$\Hilb_{m,m-1}(\tilde{R}/B)$,  as formal scheme along
$\Hilb^0_{m,m-1}(R)$, is formally smooth and $(m+1)$-dimensional
over $\k$. The natural map
$$\Hilb_{m,m-1}(\tilde{R}/B)\to \Hilb_{m-1}(\tilde{R}/B)$$
is a flat, locally complete intersection morphism of relative
dimension 1.

\endproclaim\demo{Proof} Assertion (i) is an elementary consequence
of the analysis in the proof of Theorem 1 and its proof will be
omitted. Assertion (ii) and the set-theoretic portion of Assertion
(iii) follow from this. To complete the proof of (iii),(iv), it
remains to analyze the situation locally at each pair $(I,
I')\in$Hilb$^0_{m,m-1}(R).$ We will consider the case of $(Q^m_i,
Q^{m-1}_i), 1<i<m$, as other cases are similar or simpler. There
we will focus mainly on Assertion (iv), as (iii) is essentially a
special case of this (as will be indicated below). Consider then a
pair $(I_S< I'_S)$ flatly deforming $(Q^m_i, Q^{m-1}_i)$ relative
to $B$. Then we may assume that for some $s\in\m_S, I_S$ and
$I'_S$ are generated by $xy-s$ and $f,g$ (resp. $f', g'$) with
$$f=x^{m+1-i}+\sum\limits_{j=0}^{m-i}
a_jx^j+\sum\limits_{j=1}^{i-1}b_jy^j, $$$$g=
y^i+\sum\limits_{j=0}^{m-i}
c_jx^j+\sum\limits_{j=1}^{i-1}d_jy^j,$$
$$f'=x^{m-i}+\sum\limits_{j=0}^{m-i-1}
a'_jx^j+\sum\limits_{j=1}^{i-1}b'_jy^j, $$$$g'=
y^i+\sum\limits_{j=0}^{m-i-1}
c'_jx^j+\sum\limits_{j=1}^{i-1}d'_jy^j.$$ (For the non-relative
case we take s=0.) \comment
 and in the punctual case we also set
$$a_j=d_j=a'_j=d'_j=0, \forall j,$$
$$b_j=b'_j=0, \forall j<i-1,$$
$$c_j=0, \forall j<m-i, c'_j=0, \forall j<m-i-1).$$
\endcomment

As we saw above, the relations as in (7), or equivalently (8),
are necessary and sufficient so that $I_S, I'_S$ are $S$-flat
deformations of $I, I'$ respectively. In particular, these include
$$b_{i-1}c_{m-i}=s=b'_{i-1}c'_{m-i-1}.\tag 14$$
The other relations can be used to eliminate some of the
parameters. It remains to account for the condition that
$I_S<I'_S.$ To this end it suffices to note that $$1,x, ...,
x^{m-i-1}, y,..., y^{i-1}$$ form an $S$-free basis of $R_S/I'_S$,
then express $f, g$ in terms of this basis and equate the
coefficients to 0. This yields the coefficient relations
$$a_0=-(a_{m-i}-a'_{m-i-1})a'_0+sb'_1,$$
$$a_j=a'_{j-1}-(a_{m-i}-a'_{m-i-1})a'_j, j=1,..., m-i-1,$$
$$b_j=((a_{m-i}-a'_{m-i-1})b'_j+sb'_{j+1}, j=1,...,i-2,$$
$$b_{i-1}=(a_{m-i}-a'_{m-i-1})b'_{i-1},$$
$$c_j=c'_j+c_{m-i}a'_j, j=0,..., m-i-1,$$
$$d_j=d'_j+c_{m-i}b'_j, j=1,...,i-1.$$
These coefficient relations are equivalent to
$$f=(x+a_{m-i}-a'_{m-i-1})f', g=g'+c_{m-i}f'.\tag 14.1$$
By formal manipulations, these relations
imply that
$$c'_{m-i-1}=c_{m-i}(a_{m-i}-a'_{m-i-1}).\tag 15$$
Therefore the 2 relations (14) are replaced by the single relation
$$(a_{m-i}-a'_{m-i-1})b'_{i-1}c_{m-i}=s.\tag 16$$
Consequently the relative flag Hilbert scheme is smooth here, with
regular parameters
$$a'_1,...,a'_{m-i-1}, a_{m-i}, d'_1,...,d'_{i-1}, b'_{i-1},
c_{m-i}$$ and its fibre, i.e. $\Hilb_{m,m-1}(R)$, is the
normal-crossing triple point
$$(a_{m-i}-a'_{m-i-1})b'_{i-1}c_{m-i}=0.$$
The 3 components are: $D^{m,m-1}_{i-1,i-1},$ defined by
$a_{m-i}-a'_{m-i-1}=0$  (which implies $b_{i-1}=c'_{m-i-1}=0 $);
$D^{m,m-1}_{i,i},$ defined by $b'_{i-1}=0$ (which implies
$b_{i-1}=0$); $D^{m,m-1}_{i,i-1},$ defined by $c_{m-i}=0,$ (which
implies $c'_{m-i-1}=0$).  Finally the relation (15) exhibits
$\Hilb_{m,m-1}(\tilde{R}/B)$ locally as a conic in an $\A^2$ with
coordinates $a_{m-i}, c_{m-i}$ over $\Hilb_{m-1}(\tilde{R}/B)$,
and therefore the projection is a flat locally complete
intersection morphism.\qed\enddemo\remark{Remark 5.1} In the case
of the punctual flag Hilbert  scheme, we have by Remark 2.1 that
the only coefficients in the above calculation not automatically 0
are $b_{i-1}, c_{m-i}, b'_{i-1}, c_{m-1-i}$, and the vanishing of
$a_{m-i}, a'_{m-i-1}, d_1, d'_1$ yield the relations
$$b_{i-1}=c'_{m-i-1}=c_{m-i}b'_{i-1}=0.\tag 17$$ Thus
$\Hilb^0_{m,m-1}(R)$ has two components at $(Q^m_i, Q^{m-1}_i)$:
one where $c_{m-i}=0$, which projects to
$\{Q^m_i\}\subset\Hilb_m(R)$ and to $C^{m-1}_{i-1}\subset
\Hilb_{m-1}(R);$ the other where $b'_{i-1}=0$ which projects to
$C^m_i\subset\Hilb_m(R)$ and to
$\{Q^{m-1}_i\}\subset\Hilb_{m-1}(R).$ Thus we recover in terms of
equations the set-theoretic picture presented in Theorem
5(ii).\qed\endremark \remark{Remark 5.2}For future reference we
note that in the analogous case of (relative) deformations of
$(Q^m_i, Q^{m-1}_{i-1})$, $f', g'$ take the form
$$f'=x^{m-i+1}+\sum\limits_{j=0}^{m-i}
a'_jx^j+\sum\limits_{j=1}^{i-2}b'_jy^j, g'=
y^{i-1}+\sum\limits_{j=0}^{m-i}
c'_jx^j+\sum\limits_{j=1}^{i-2}d'_jy^j.$$ The relations (15) and
$b_{i-1}=(a_{m-i}-a'_{m-i-1})b'_{i-1}$ (see just above (15)) are
replaced by
$$c_{m-i}=(d_{i-1}-d'_{i-2})c'_{m-i},
b'_{i-2}=(d_{i-1}-d'_{i-2})b_{i-1}.\tag 18$$ As we shall see, this
remark is useful in studying flag Hilbert schemes with more than 2
constituents.\qed\endremark

\remark{Remark 5.3} Note that for any $(I,
I')\in\Hilb_{m,m-1}(R)$, the annihilator ${\text {Ann}}(I'/I)<R$
is an ideal of colength 1. This yields a morphism
$$A_m:\Hilb_{m,m-1}(R)\to\Hilb_1(R)=\Spec(R).$$ For example, in the
situation of (?), we have by (the analogue of) (8),
$$yf'=b'_{i-1}g'=b'_{i-1}g-b'_{i-1}c_{m-i}f',$$
hence $J:=(x+a_{m-i}-a'_{m-i-1}, y+b'_{i-1}c_{m-i})$ annihilates
$f' \mod I_S.$  Hence by (?) again, $J$ annihilates $I'_S \mod
I_S.$ Note that by (?), we have $b'_{i-1}c_{m-i}=
d_{i-1}-d'_{i-1}.$ Therefore the value of $A_m$ on the $S-$valued
point $(I_S, I'_S)$ (which extends the closed point $(Q^m_i,
Q^{m-1}_i)$) is
$$A_m(I_S, I'_S)=(x+a_{m-i}-a'_{m-i-1}, y+d_{i-1}-d'_{i-1}).$$
Consequently, the closed (special) fibre of $A_m$ is defined in
terms of $f,g,f',g'$ by the condition that
$$f(x,0)=xf'(x,0), g(0,y)=g'(0,y).\qed$$

\endremark
\remark{Construction 5.4} An analogue of Construction 3.2 in the
flag case may be given as follows. Let $a_i, d_i$ etc be as there
and let $a'_i, d'_i$ etc. be the analogous objects for
$\Hilb_{m-1}$. Set
$$r=a_{m-1}-a'_{m-2}, s=d_{m-1}-d'_{m-2},$$ so that
$$F_0=F'_0(x+r), F_m=F'_{m-1}(y+s)$$
and we have the relation $$rs=t,$$ so that $(r,s)$ give a copy of
$\Spec(\tilde{R})/B$.  Then $\Hilb_{m,m-1}(\tilde{R}/B)$ may be
relalized as the subscheme of
$\tilde{C}\times_B\Spec(\tilde{R})\times_B\Hilb_{m-1}(\tilde{R}/B)$
defined by $$u'_iv_i=ru_iv'_i,\ \  v'_iu_{i+1}=sv_{i+1}u'_i$$
$$a'_iu_{i+1}=sd'_{m-1-i}v_{i+1},\ \  d'_{m-1-i}v_i=ra'_iu_i,
i=1,...,m-1.$$

\qed
\endremark

 In extending these results to the case of
longer- a fortiori, full- flags, the same methods apply. But in
the conclusions, a couple of new twists come up, already for flags
of type $m.=(m,m-1,m-2).$ First, the punctual Hilbert scheme
$\Hilb^0_{m.}(R)$ will have components of varying dimensions (in
this case, 1 and 2), roughly speaking because the parameters in
$\Hilb^0_m(R)$ and $\Hilb^0_{m-2}(R)$ can vary independently.
Second, the relative flag Hilbert scheme $\Hilb_{m.}(\tilde{R}/B)$
is not the same locally at $(Q^m_i, Q^{m-1}_i, Q^{m-2}_i)$ as at
$(Q^m_i, Q^{m-1}_i, Q^{m-2}_{i-1})$. The following proof contains
the relevant computation. \proclaim{ Lemma 6} Set
$m.=(m,m-1,m-2).$ Then (i) as algebraic set, $\Hilb^0_{m.}(R)$ is
of the form
$$ C^m_1\cup C^{m-1}_1\cup C^{m,m-2}_{2,1}\cup C^{m-1}_2\cup
...\cup C^{m,m-2}_{m-2, m-3}\cup C^{m-1}_{m-2}\cup C^m_{m-1}.$$
Each component $C^{m,m-2}_{i, i-1}$ projects isomorphically to
$C^m_i\times C^{m-2}_{i-1}\subset\Hilb_{m,m-2}(R)$ and to
$\{Q^{m-1}_i\}\subset\Hilb_{m-1}(R)$.\par (ii)
$\Hilb_{m.}(\tilde{R}/B)$ is irreducible and is smooth except at
points $(Q^m_i, Q^{m-1}_i, Q^{m-2}_{i-1})$, where is has a rank-4
quadratic hypersurface singularity with local equation
$$
(a_{m-i}-a'_{m-i-1})c_{m-i}=(d'_{i-1}-d''_{i-2})c''_{m-i-1}\tag
19$$

\endproclaim
\demo{proof} The set-theoretic assertion (i) follows from Theorem
5. For (ii), we will analyze $\Hilb_{m.}(\tilde{R}/B)$ locally at
points of the form $(Q^m_i, Q^{m-1}_i, Q^{m-2}_{i})$ or $(Q^m_i,
Q^{m-1}_i, Q^{m-2}_{i-1})$ as other cases are similar or simpler.
Beginning with the former case, consider a deformation
$(I''_S<I'_S<I_S)$ of $(Q^m_i, Q^{m-1}_i, Q^{m-2}_{i})$ where
$I_S, I'_S$ are as in the proof of Theorem 5 and $I''_S$ is
analogously defined by
$$f''=x^{m-i-1}+\sum\limits_{j=0}^{m-i-2}a''_jx^j+
\sum\limits_{j=1}^{i-1}b''_jy^j,$$
$$g''=y^{i}+\sum\limits_{j=0}^{m-i-2}c''_jx^j+
\sum\limits_{j=1}^{i-1}d''_jy^j.$$ Then working as in the proof of
Theorem 5 we find the the $(m+1)$ parameters
$$a''_1,...,a''_{m-i-2}, d''_1,...,d''_{i-1}, b''_{i-1}, c_{m-i},
a'_{m-i-1}, a_{m-i}$$ such that all the coefficients of all our
polynomials $f,...,g''$, as well as the parameter $s$, are regular
expressions in these and there are no relations. This shows that
$\Hilb_{m.}(\tilde{R}/B)$ is smooth at $(Q^m_i, Q^{m-1}_i,
Q^{m-2}_{i})$.\par In the case of $(Q^m_i, Q^{m-1}_i,
Q^{m-2}_{i-1})$,  we may assume
$$f''=x^{m-i}+\sum\limits_{j=0}^{m-i-1}a''_jx^j+
\sum\limits_{j=1}^{i-2}b''_jy^j,$$
$$g''=y^{i-1}+\sum\limits_{j=0}^{m-i-1}c''_jx^j+
\sum\limits_{j=1}^{i-2}d''_jy^j$$ and analogous considerations
yield $(m+2)$ parameters
$$a''_1,...,a''_{m-i-1}, d''_1,...,d''_{i-2}, b'_{i-1}, c_{m-i},
c''_{m-i-1}, d'_{i-1}, a_{m-i}$$ such that all the coefficients of
all our polynomials, as well as the parameter $s$, are regular
expressions in these, and satisfying the relation
$$c'_{m-i-1}= (a_{m-i}-a'_{m-i-1})c_{m-i}=
(d'_{i-1}-d''_{i-2})c''_{m-i-1}$$ which yields the equation
(19).\par \comment Finally, since
$$s=b'_{i-1}(a_{m-i}-a'_{m-i-1})c_{m-i},$$ it is a product of 3
distinct members of a system of parameters, hence a nonzero
divisor, hence $\Hilb_{m.}(\tilde{R}/B)$ has, locally at $(Q^m_i,
Q^{m-1}_i, Q^{m-2}_{i-1})$, no component not dominating $B$ and
since $\Hilb_{m.}(\tilde{R}/B)$ is clearly generically irreducible
over $B$, it follows that it is irreducible.
\endcomment
Finally, irreducibility of $\Hilb_{m.}(\tilde{R}/B)$ follows from
the fact that its natural map to $\Hilb_{m,m-1}(\tilde{R}/B)$ is
flat with irreducible generic fibre and
$\Hilb_{m,m-1}(\tilde{R}/B)$ is irreducible by Theorem 5(iv).\qed
\enddemo
 \proclaim{Theorem 7} The (full) flag
Hilbert scheme $fHilb_m(\tilde{R}/B)$ has locally complete
intersection singularities and its natural map to $B$ is a local
complete intersection morphism. In particular
$fHilb_m(\tilde{R}/B)$ is reduced and is flat over $B$.
\endproclaim
\demo{proof} The fact that $\fHilb_m(\tilde{R}/B)$ has locally
complete intersection singularities follows by a downwards
induction as in the proof of Lemma 6. A similar argument shows
that each map $$\fHilb_m(\tilde{R}/B)\to
\fHilb_{m-1}(\tilde{R}/B)$$ is a locally complete intersection
morphism, hence so is $\fHilb_m(\tilde{R}/B)\to B.$\qed\enddemo As
in Remark 5.3, we have a natural map
$$fA_m=A_m\times A_{m-1}\times...\times A_1:\fHilb_m(R)\to\Hilb_1(R)^m=\Spec(R)^m$$
which fits in the diagram $$\matrix \fHilb_m(R)&\to & \Hilb_m(R)\\
fA_m\downarrow&&\downarrow {\text{cyc}}\\ \Spec(R)^m&\to &
\Sym^m(\Spec(R))\endmatrix$$ Then the closed fibre of $fA_m$ is
called the $m$th {\it{punctual flag-Hilbert scheme}} of $R$,
denoted $\fHilb^0_m(R)$, and the (scheme-theoretic) projection of
$\fHilb^0_m(R)$ to $\Hilb_m(R)$, i.e. the closed fibre of the
cycle map cyc, is called the $m$th {\it{punctual Hilbert scheme}}
of $R$, denoted $\Hilb^0_m(R).$ It is clear that $\fHilb^0_m(R)$
and $\Hilb^0_m(R)$ endow the similarly denoted algebraic sets
discussed previously with a scheme structure.
 \proclaim{Corollary 8} The
full flag punctual Hilbert scheme
$\fHilb^0_m(R)=\Hilb^0_{m,...,1}(R)$ is reduced and is the
transverse union of components of the form
$$C^{m,m-2,...,m-2j}_{i, i-1, ...,i-j}, \forall i, 1\leq i\leq
m-1, j=\min([\frac{m}{2}], i-1, m-i-1)\geq 0,$$ which projects
isomorphically to $C^m_i\times ...\times C^{m-2j}_{i-j}$, and to a
point in the other factors;
$$C^{m-1,m-2,...,m-2j-1}_{i, i-1, ...,i-j}, \forall i, 1\leq i\leq m-2,
j=\min([\frac{m-1}{2}], i-1, m-i-1-2)\geq 0,$$ which projects
isomorphically to $C^{m-1}_i\times ...\times C^{m-2j-1}_{i-j}$ and
to a point in the other factors. The punctual Hilbert scheme
$Hilb^0_m(R)$, with the scheme structure as above is a reduced
nodal curve.
\endproclaim\demo{proof} Except for the assertion that  $\fHilb^0_m(R)$,
hence also $\Hilb^0_m(R)$ are reduced, this is a straightforward
extension of Lemma 6; the constraints on $j$ come simply from the
fact that the components of $\Hilb^0_k(R)$ are $C^k_j,
j=1,...,k-1$. For the reducedness assertion, we argue by induction
on $m$. We may work locally at a point of the form
$(Q^m_i,Q^{m-1}_i,...)$, as other cases are similar or simpler.
Consider a deformation of the form $(I_S=(f,g)< I'_S=(f',g')<...)$
with $I_S, I'_S$ as in ??, that yields an $S$-valued point of
$\fHilb^0_m(R)$. By induction, we may assume $$f'(x,0)=x^{m-i},
g'(0,y)=y^i,$$ which implies that $$f'=x^{m-i}+b'_{i-1}y^{i-1},
g'=y^i+c'_{m-i-1}x^{m-i-1}.$$ Then by Remark 5.3 we get a similar
conclusion for $f,g.$ Thus , in terms of the coordinates
$a_1,...,a_{m-i}, d_1,...,d_{i-1}, b_{i-1}, c_{m-i}$ as in the
proof of Theorem 2, the {\it subscheme} $\Hilb^0_m(R)\subset
\Hilb_m(R)$ is simply defined by the vanishing of
$a_1,...,d_{i-1}$, hence is reduced as claimed, and likewise for
$\fHilb^0_m(R)$. \qed\enddemo \remark{Remark 8.1} To 'picure' the
configuration in Corollary 7, it is amusing to display the
components of $\Hilb^0_k(R), k=2,...,m $ as segments arranged as
an isosceles triangle:$$-
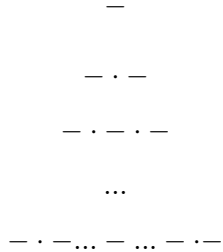
$$
$$-\cdot-$$$$-\cdot-\cdot-$$$$...$$$$-\cdot-...-...-\cdot-$$
 Then the components of
$\fHilb^0_m(R)$ are the columnwise products of these components.

\endremark
\comment
 \proclaim{Theorem 8} The (full) flag
Hilbert scheme $fHilb_m(\tilde{R}/B)$ has locally complete
intersection singularities and its natural map to $B$ is a local
complete intersection morphism. In particular
$fHilb_m(\tilde{R}/B)$ is reduced and is flat over $B$.

\endproclaim
\demo{proof} The fact that $fHilb_m(\tilde{R}/B)$ has locally
complete intersection singularities follows by a downwards
induction as in the proof of Lemma 7. A similar argument shows
that each map $$\fHilb_m(\tilde{R}/B)\to
\fHilb_{m-1}(\tilde{R}/B)$$ is a locally complete intersection
morphism, hence so is $\fHilb_m(\tilde{R}/B)\to B.$\qed\enddemo
\endcomment
 \remark{Remark 9} in \cite{R} we construct, based on
'geometric' considerations, a space $W^m(X/B)$ together with a
morphism $J_m: W^m(X/B)\to \fHilb_m(X/B)$, which we will prove is
an isomorphism. This proof requires that we know a priori that
$\fHilb_m(X/B)$ is {\it reduced}. \endremark\remark{Remark 10} It
seems likely that the above results go through without the
assumption that the base field $\k$ is algebraically closed,
provided the fibre nodes of $X/B$ are of 'split' type, i.e. each
node $p$, and each of the 2 tangent directions at $p$, are defined
over $\k$. Beyond this however, it seems some further analysis is
needed. For instance, if $p$ is a node defined over $\k$ of
nonsplit type, i.e. with equation formally equivalent to $x^2+y^2,
\text{char}(\k)\neq 2,$ the punctual Hilbert scheme $\Hilb^0_m$ at
$p$ apparently has just 1 or 0 components defined over $\k$
depending on whether $m$ is even or odd; the other components
occur in conjugate pairs.\endremark

\enddocument